\documentclass[12pt]{amsart}
\usepackage{amssymb,amsmath,amsthm}
\usepackage{tikz}
\usepackage{cancel}
\usepackage{mathabx}

\usetikzlibrary{arrows,shapes,automata,backgrounds,decorations,petri,positioning}

\theoremstyle{plain}
\newtheorem{thm}{Theorem}[section]
\newtheorem{lem}[thm]{Lemma}

\newtheorem{prop}[thm]{Proposition}
\theoremstyle{definition}
\newtheorem{defn}[thm]{Definition}

\newtheorem{ex}[thm]{Example}

\newcommand{\mf}[1]{\mbox{$\mathfrak #1$}}
\newcommand{\mc}[1]{\mbox{$\mathcal #1$}}
\newcommand{\ts}[1]{\mbox{$\textsf{#1}$}}

\renewcommand{\star}{\mbox{$\Asterisk$}}

\title{A combinatorial proof of symmetry among minimal star factorizations}

\author{Bridget Eileen Tenner}
\address{Department of Mathematical Sciences, DePaul University, Chicago, IL 60614}
\email{bridget@math.depaul.edu}


\begin{document}

\begin{abstract}
The number of minimal transitive star factorizations of a permutation was shown by Irving and Rattan to depend only on the conjugacy class of the permutation, a surprising result given that the pivot plays a very particular role in such factorizations.  Here, we explain this symmetry and provide a bijection between minimal transitive star factorizations of a permutation $\pi$ having pivot $k$ and those having pivot~$k'$.\\

\end{abstract}

\maketitle

\section{Introduction}

For any positive integer $n$, let $[n]$ denote the set $\{1, \ldots, n\}$.  The symmetric group $\mf{S}_n$ consists of all permutations of the set $[n]$.  It is most convenient for our purposes to represent permutations in cycle notation.  That is, write a permutation $\pi$ as the product of disjoint cycles $(x, \pi(x), \ldots, \pi^{\ell-1}(x))$ where $\ell>0$ is minimal so that $\pi^{\ell}(x) = x$.  If an element $x$ is fixed by $\pi$, that is, if $\pi(x) = x$, then the cycle $(x)$ may  be suppressed.

\begin{ex}\label{ex:cycle notation}
Let $\pi \in \mf{S}_6$ be the permutation defined as follows. 
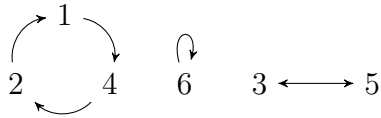
\begin{center}
\begin{tikzpicture}[node distance=.5cm,>=stealth',bend angle=45,auto]
\tikzstyle{state}=[minimum size=4mm]
\node[state] (a) {$1$}; \node[state] [below right=of a, xshift=-2.5mm] (b) {$4$}; \node[state] [left=of b, xshift=-2.5mm] (c) {$2$};
\path[->] (a) edge[bend left] (b); \path[->] (b) edge[bend left] (c); \path[->] (c) edge[bend left] (a);
\node[state] (e) [right=of b] {$6$};
\path[->] (e) edge [loop above] (e);
\node[state] (f) [right=of e] {$3$}; \node[state] (g) [right=of f, xshift=5mm] {$5$};
\path[<->] (f) edge (g);
\end{tikzpicture}
\end{center}
Then $\pi = (421)(6)(35) = (6)(214)(53)$.
\end{ex}

As evidenced by Example~\ref{ex:cycle notation}, there are multiple ways to write a given permutation in cycle notation.  It is common to designate one of these to be standard.

\begin{defn}
The \emph{standard form} of a permutation in cycle notation is obtained by writing each cycle so that its minimal element appears in the leftmost position, and writing the cycles from left to right in increasing order of minimal elements.  Fixed points are not suppressed in the standard form.
\end{defn}

\begin{ex}
The standard form of the permutation in Example~\ref{ex:cycle notation} is $(142)(35)(6)$.
\end{ex}

There are several well studied sets of generators for the symmetric group, and effort has been made to enumerate the most efficient (that is, the shortest) ways to write a given permutation in terms of generators in one of these sets.

\begin{defn}
Given a generating set $G$ of the symmetric group $\mf{S}_n$, a factorization, or decomposition of $\pi \in \mf{S}_n$ as $\pi = g_1 \cdots g_{\ell}$ for $g_i \in G$ is \emph{minimal} if $\ell$ is minimal.  This $\ell$ is the \emph{$G$-length} (or \emph{length} if the generating set is clear from the context) of~$\pi$.
\end{defn}

When $G = \{(i \ i+1) : 1 \le i < n\}$, Stanley has shown that in certain cases, including the $G$-longest permutation $(1n)(2(n-1))(3(n-2))\cdots$, this number of minimal factorizations of a permutation $\pi$ is the same as the number of standard Young tableaux of a particular shape $\lambda(\pi)$ \cite{stanley}.  In a different vein, D\'enes computed the number of minimal factorizations of a permutation when $G = \{(i \ j) : 1 \le i<j \le n\}$~\cite{denes}.

The generating set for $\mf{S}_n$ that we consider here is the set of star transpositions.

\begin{defn}\label{defn:star}
Fix positive integers $n \ge k \ge 1$.  The set $\star_{n;k} = \{(k \ i) : i \in [n] \setminus \{k\}\}$ is the set of \emph{star transpositions with pivot~$k$}.
\end{defn}

The nomenclature refers to the fact that when the elements of $\star_{n;k}$ are considered to define edges on the vertices $[n]$, the resulting graph is a star with center label~$k$.

\begin{lem}
The set $\star_{n;k}$ generates the symmetric group $\mf{S}_n$, for any $k \in [n]$.
\end{lem}

\begin{proof}
Since $(i \ i+1) = (k \ i)(k \ (i+1))(k \ i)$, we can generate all simple transpositions, and these, in turn, generate~$\mf{S}_n$.
\end{proof}

In \cite{pak}, Pak considered minimal factorizations of a particular family of permutations into the star transpositions of Definition~\ref{defn:star}.

\begin{defn}
A factorization $\pi = g_1 \cdots g_r \in \mf{S}_n$ for $g_i \in \star_{n;k}$ is \emph{transitive} if the group generated by $\{g_1, \ldots, g_r\}$ acts transitively on the set $[n]$.  In other words, this factorization is transitive if $\{g_1, \ldots, g_r\} = \star_{n;k}$.
\end{defn}

In the class of permutations studied in \cite{pak}, the only fixed point was the pivot itself.  Thus the factorizations of these permutations are necessarily transitive.  Pak's work was generalized by Irving and Rattan, who computed the number of minimal transitive star factorizations of any permutation~\cite{irving-rattan}.  In that work, Irving and Rattan discovered a surprising symmetry in their enumeration, essentially saying that the choice of pivot does not affect the number of minimal transitive star factorizations of a permutation.  The purpose of the current article is to provide a combinatorial proof of the symmetry that they found.

\begin{defn}
Given a permutation $\pi \in \mf{S}_n$, let $\star_k(\pi)$ denote the set of minimal transitive star factorizations of $\pi$ having pivot $k$.  Let $s_{k}(\pi) = |\star_{k}(\pi)|$.
\end{defn}

As discussed in \cite{irving-rattan}, if $\pi \in \mf{S}_n$ has $m$ cycles, then each element of $\star_k(\pi)$ has length $n+m-2$.  As is customary, permutations are viewed as maps, and so are multiplied from right to left.

\begin{ex}
Consider $(142)(35)(6) \in \mf{S}_6$.  Then
$$(31)(36)(36)(32)(34)(31)(35) \in \star_{3}\big((142)(35)(6)\big).$$
\end{ex}

The main result of \cite{irving-rattan} is the following theorem, where $x_{(r)}$ denotes the \emph{falling factorial}
$$x_{(r)} = x(x-1)(x-2)\cdots(x-r+1) = \frac{x!}{(x-r)!}.$$
Note that the discussions in \cite{irving-rattan} assume that $k=1$, but the proof can be extended to an arbitrary pivot~$k$.

\begin{thm}[\cite{irving-rattan}]\label{thm:irving-rattan}
Let $\pi \in \mf{S}_n$ have cycles of lengths $\ell_1, \ldots, \ell_m$.  Then
$$s_{k}(\pi) = (n+m-2)_{(m-2)}\ell_1\cdots\ell_m$$
for all $k \in [n]$.
\end{thm}

The remarkable feature of Theorem~\ref{thm:irving-rattan} is its symmetry: it depends only on the cycle type of $\pi$.  Given the special role played by the pivot in star transpositions, one would not expect the cycle containing the pivot to behave in the same way as the other cycles in $\pi$.  It was an open question in \cite{irving-rattan} to explain this symmetry.  Goulden and Jackson looked at non-minimal factorizations into star transpositions, and uncovered the same symmetry in that setting \cite{goulden-jackson}, although again an explanation of this symmetry was lacking.  In \cite{feray}, F\'eray has given a proof of the symmetry in both the minimal and non-minimal situations.  His argument uses the algebra of partial permutations of Ivanov and Kerov \cite{ivanov-kerov}, but does not give a combinatorial reason for the phenomenon.

Giving such a justification, in fact a bijection between minimal transitive star factorizations with pivot $k$ and those with pivot $k'$, is the purpose of the current article.  We do this by giving a bijection $\phi_{\pi,k}$ between injections from $[m-2]$ into $[n+m-2]$ together with elements of $[\ell_1] \times \cdots \times [\ell_m]$, and the elements of $\star_k(\pi)$ (see Definition~\ref{defn:phi} and Theorem~\ref{thm:bijection}).  This $\phi_{\pi,k}$ takes such an injection and $m$-tuple, and produces a \emph{valid word} for $\star_k(\pi)$ via the maps $\ts{tree}_{\pi,k}$ and $\omega$, and a set of \emph{cycle enclosures} via the map $\ts{cycle}_{\pi,k}$.  As shown in Proposition~\ref{prop:enough information}, these determine a unique element of $\star_k(\pi)$.

We can compose the maps
$$\phi_{\pi,k'} \circ \left(\phi_{\pi,k}\right)^{-1}$$
to obtain the desired combinatorial bijection between $\star_k(\pi)$ and $\star_{k'}(\pi)$.  We will often describe the inverses of our procedures, such as in Proposition~\ref{prop:tree is a bijection}, to give intuition about how to apply the above composition of maps.

In Section~\ref{section:notation}, we briefly establish notation that will be used in the duration of the paper.  Section~\ref{section:characterization} characterizes elements of $\star_k(\pi)$, relying heavily on the work of \cite{irving-rattan}.  It is also in this section that we define cycle enclosures and valid words.  In Section~\ref{section:trees}, we introduce the class of trees which are crucial to our bijection, and which themselves are in bijection with the set of valid words via a map $\omega$.  In Section~\ref{section:construction}, we define the maps $\ts{cycle}_{\pi,k}$ and $\ts{tree}_{\pi,k}$, and ultimately the bijection $\phi_{\pi,k}$.  Finally, the paper concludes with Section~\ref{section:conclusion}, in which we give intuition for understanding the symmetry of Theorem~\ref{thm:irving-rattan}.  Throughout the paper, we use running examples to illustrate each of the definitions and operations.

\section{Notation and terminology}\label{section:notation}

Here we establish  notation and terminology that will be used throughout the present work.

Fix positive integers $n \ge k$, and $\pi \in \mf{S}_n$.  Let $\pi$ consist of $m$ disjoint cycles, with lengths $\ell_1, \ldots, \ell_m$ when read from left to right in standard form.  The symbol $k$ appears in the $p$th of these cycles.  (Note that $k$ will play the role of the pivot, hence the index ``$p$.'')

Each tree discussed in this work is \emph{ordered}: it has a designated root node, and an ordering is specified for the children of each vertex.

\section{Characterization of minimal transitive star factorizations}\label{section:characterization}

In this section we  give a description of $\star_k(\pi)$.  The interested reader is referred to \cite{irving-rattan} for more information.  While the discussions in \cite{irving-rattan} assume that $k=1$, the proofs of the results cited below can be extended to an arbitrary pivot value~$k$.

\begin{ex}\label{ex:mtsf ex}
Consider the permutation $\pi = (142)(35)(6) \in \mf{S}_6$.  By Theorem~\ref{thm:irving-rattan},
$$s_{1}(\pi) = \frac{(6+3-2)!}{6!}\cdot 3\cdot 2 \cdot 1 = 42.$$
The elements of $\star_{3}(\pi)$, that is, the $42$ minimal transitive star factorizations with pivot $3$ of $\pi$,
are given below.
$$\begin{array}{ccc} 
(35)(36)(36)(31)(32)(34)(31) & (35)(36)(36)(32)(34)(31)(32) & (35)(36)(36)(34)(31)(32)(34)\\
(36)(36)(35)(31)(32)(34)(31) & (36)(36)(35)(32)(34)(31)(32) & (36)(36)(35)(34)(31)(32)(34)\\
(36)(36)(31)(32)(34)(31)(35) & (36)(36)(32)(34)(31)(32)(35) & (36)(36)(34)(31)(32)(34)(35)\\
(35)(36)(31)(32)(34)(31)(36) & (35)(36)(32)(34)(31)(32)(36) & (35)(36)(34)(31)(32)(34)(36)\\
(36)(31)(32)(34)(31)(36)(35) & (36)(32)(34)(31)(32)(36)(35) & (36)(34)(31)(32)(34)(36)(35)\\
(35)(31)(36)(36)(32)(34)(31) & (35)(32)(36)(36)(34)(31)(32) & (35)(34)(36)(36)(31)(32)(34)\\
(31)(36)(36)(32)(34)(31)(35) & (32)(36)(36)(34)(31)(32)(35) & (34)(36)(36)(31)(32)(34)(35)\\
(35)(31)(32)(36)(36)(34)(31) & (35)(32)(34)(36)(36)(31)(32) & (35)(34)(31)(36)(36)(32)(34)\\
(31)(32)(36)(36)(34)(31)(35) & (32)(34)(36)(36)(31)(32)(35) & (34)(31)(36)(36)(32)(34)(35)\\
(35)(31)(32)(34)(36)(36)(31) & (35)(32)(34)(31)(36)(36)(32) & (35)(34)(31)(32)(36)(36)(34)\\
(31)(32)(34)(36)(36)(31)(35) & (32)(34)(31)(36)(36)(32)(35) & (34)(31)(32)(36)(36)(34)(35)\\
(35)(31)(32)(34)(31)(36)(36) & (35)(32)(34)(31)(32)(36)(36) & (35)(34)(31)(32)(34)(36)(36)\\
(31)(32)(34)(31)(35)(36)(36) & (32)(34)(31)(32)(35)(36)(36) & (34)(31)(32)(34)(35)(36)(36)\\
(31)(32)(34)(31)(36)(36)(35) & (32)(34)(31)(32)(36)(36)(35) & (34)(31)(32)(34)(36)(36)(35)
\end{array}$$
Note that, in some of these, the identity product $(36)(36) = (1)$ is included.  This is done so that the entire product factorization is transitive.
\end{ex}

The following statements are easy to prove, and are discussed in \cite{irving-rattan}.

\begin{lem}[\cite{irving-rattan}]\label{lem:allowable cycle decomps}
\begin{enumerate}\renewcommand{\labelenumi}{(\alph{enumi})}
\item The cycle
$$(k \ a_2 \ a_3 \ \cdots \ a_{\ell})$$
admits exactly one minimal $k$-star factorization:
$$(k \ a_{\ell})(k \ a_{\ell-1}) \cdots (k \ a_3)(k \ a_2).$$
\item The cycle
$$(b_1 \ b_2 \ \cdots \ b_{\ell}),$$
where $b_i \neq k$ for all $i$, admits $\ell$ different minimal $k$-star factorizations:
$$(k \ b_i)(k \ b_{i+\ell-1})(k \ b_{i+\ell-2}) \cdots (k \ b_{i+1})(k \ b_i),$$
where the subscripts are taken modulo~$\ell$.
\end{enumerate}
\end{lem}

It is helpful to introduce terminology to identify the different possibilities described in Lemma~\ref{lem:allowable cycle decomps}(b).

\begin{defn}
Suppose that the standard form of a permutation $\pi$ contains the cycle $C = (b_1 \ b_2 \ \cdots \ b_{\ell})$, and suppose that $\delta \in \star_{k}(\pi)$, with $k \neq b_j$ for all $j$.  Let $i$ be such that $\delta$ contains the subword
$$(k \ b_i)(k \ b_{i+\ell - 1}) \cdots (k \ b_{i+1})(k \ b_i),$$
with subscripts taken modulo $\ell$.  (Note that the transpositions in this subword do not necessarily appear consecutively in the factorization $\delta$, as shown in Example~\ref{ex:cycle enc}.)  Then the cycle $C$ is \emph{enclosed} by~$b_i$.
\end{defn}

\begin{defn}
Given a permutation $\pi$ and ${\delta} \in \star_{k}(\pi)$, the \emph{cycle enclosures} of ${\delta}$ are the set of letters that enclose all cycles in $\pi$ except the cycle containing~$k$.
\end{defn}

\begin{ex}\label{ex:cycle enc}
Continuing our running example, consider ${\delta} = (35)(34)(31)(32)(36)(36)(34) \in \star_{3}((142)(35)(6))$.  Then $(142)$ is enclosed by $4$, and the cycle enclosures of ${\delta}$ are $\{4,6\}$. 
\end{ex}

Before stating Lemma~\ref{lem:allowable subwords}, which is crucial to the description of $\star_k(\pi)$, we must make the following definition.

\begin{defn}\label{defn:decomp to word}
Given ${\delta} = (k \ {\delta}_1)(k \ {\delta}_2) \cdots (k \ {\delta}_{n+m-2}) \in \star_{k}(\pi)$, define a word $\omega({\delta}) \in [m]^{n+m-2}$ so that if ${\delta}_i$ appears in the $j$th cycle in the standard form of $\pi$, then the $i$th letter of $\omega(\delta)$ is~$j$.
\end{defn}

\begin{ex}\label{ex:delta(omega)}
Rewrite each ${\delta} \in \star_{3}((142)(35)(6))$ from Example~\ref{ex:mtsf ex} as the word $\omega({\delta})$.
$$\begin{array}{ccc} 
2331111 & 2331111 & 2331111\\
3321111 & 3321111 & 3321111\\
3311112 & 3311112 & 3311112\\
2311113 & 2311113 & 2311113\\
3111132 & 3111132 & 3111132\\
2133111 & 2133111 & 2133111\\
1331112 & 1331112 & 1331112
\end{array}
\hspace{1in}
\begin{array}{ccc}
2113311 & 2113311 & 2113311\\
1133112 & 1133112 & 1133112\\
2111331 & 2111331 & 2111331\\
1113312 & 1113312 & 1113312\\
2111133 & 2111133 & 2111133\\
1111233 & 1111233 & 1111233\\
1111332 & 1111332 & 1111332
\end{array}$$
The repetition in this list is due to Lemma~\ref{lem:allowable cycle decomps}(b).
\end{ex}

The following lemma completely characterizes the possible words $\omega({\delta})$ that may exist for ${\delta}$ a minimal transitive star factorization of a permutation~$\pi$.

\begin{lem}[\cite{irving-rattan}]\label{lem:allowable subwords}
Let $\omega \in [m]^{n+m-2}$ be a word on $[m]$.  There exists ${\delta} \in \star_{k}(\pi)$ such that $\omega = \omega({\delta})$ if and only if the following statements hold for~$\omega$:
\begin{itemize}
\item the symbol $p$ appears $\ell_p - 1$ times,
\item the symbol $j$ appears $\ell_j+1$ times for all $j \in [m] \setminus \{p\}$,
\item the word $\omega({\delta})$ contains no subword $ijij$ for $i \neq j$, and
\item the word $\omega({\delta})$ contains no subword $jpj$ for $j \neq p$.
\end{itemize}
\end{lem}

\begin{defn}
If $\omega \in [m]^{n+m-2}$ satisfies the requirements of Lemma~\ref{lem:allowable subwords},  that $\omega$ is a \emph{valid word} for $\star_{k}(\pi)$.  Let $\mc{W}_{k}(\pi)$ be the set of valid words for~$\star_{k}(\pi)$.
\end{defn}

In fact, the information contained in a minimal transitive star factorization with pivot $k$ is equivalent to the information contained in its cycle enclosures and its image under~$\omega$.  We now explain this precisely.  Encoding a minimal transitive star factorization in this way will be key to the description of the bijection~$\phi_{\pi,k}$.

\begin{defn}\label{defn:word and enclosures to mtsf}
Let $(\ell'_1, \ldots, \ell'_{m-1}) = (\ell_1, \ldots, \widehat{\ell_p},\ldots,\ell_m)$.  Define a map
$$\rho_{\pi,k} \ : \ \mc{W}_{k}(\pi) \times [\ell'_1] \times \cdots \times [\ell'_{m-1}] \ \rightarrow \ \star_{k}(\pi)$$
as follows.  Consider $(\omega,c_1,\ldots,c_{m-1}) \in \mc{W}_{k}(\pi) \times [\ell'_1] \times \cdots \times [\ell'_{m-1}]$.  Write the $p$th cycle in the standard form of $\pi$ as $(k \ a_2 \ \cdots \ a_{\ell_p})$.  Replace the $\ell_p-1$ copies of $p$ in $\omega$ by the star transpositions
$$(k \ a_{\ell_p})(k \ a_{\ell_p-1}) \cdots (k \ a_3)(k \ a_2),$$
in order.

For $i<p$, suppose that the $i$th cycle of $\pi$ when written in standard form is $(b_1 \ b_2 \ \cdots \ b_{\ell_i})$.  Replace the $\ell_i'+1 = \ell_i+1$ copies of $i$ in $\omega$ by the star transpositions
$$(k \ b_{c_i})(k \ b_{c_i +\ell_i - 1}) \cdots (k \ b_{c_i+1})(k \ b_{c_i}),$$
in order, where the subscripts are taken modulo~$\ell_i$.

For $i>p$, suppose that the $i$th cycle of $\pi$ when written in standard form is $(b_1 \ b_2 \ \cdots \ b_{\ell_i})$.  Replace the $\ell'_{i-1}+1 = \ell_i+1$ copies of $i$ in $\omega$ by the star transpositions
$$(k \ b_{c_{i-1}})(k \ b_{c_{i-1} + \ell_i - 1}) \cdots (k \ b_{c_{i-1}+1})(k \ b_{c_{i-1}}),$$
in order, where the subscripts are taken modulo~$\ell_i$.

This uniquely determines a minimal transitive star factorization with pivot $k$, which we denote $\rho_{\pi,k}(\omega,c_1, \ldots, c_{m-1})$.
\end{defn}

\begin{ex}
Take $2111331 \in \mc{W}_{3}((142)(35)(6))$, $2 \in [3]$, and $1 \in [1]$.  Then
$$\rho_{(142)(35)(6),3}(2111331,2,1) = (35)(34)(31)(32)(36)(36)(34).$$
\end{ex}

\begin{prop}\label{prop:enough information}
The map $\rho_{\pi,k}$ is a bijection.
\end{prop}

\begin{proof}
The operation of $\rho_{\pi,k}$ is easily reversible: take ${\delta} \in \star_{k}(\pi)$, let $\omega = \omega({\delta})$, and define the $c_i$ from the indices of the set of cycle enclosures of~$\delta$.
\end{proof}

It is clear from each of Lemmas~\ref{lem:allowable cycle decomps} and~\ref{lem:allowable subwords} that the cycle containing the pivot in a permutation $\pi$ behaves differently with regard to elements of $\star_{k}(\pi)$.  This emphasizes the unexpected nature of the symmetry in Theorem~\ref{thm:irving-rattan}.

\section{A class of trees}\label{section:trees}

In \cite{irving-rattan}, a correspondence was given between $\star_1(\pi)$ and a particular class of trees.  We will similarly utilize a graphical approach to explain the symmetry of Theorem~\ref{thm:irving-rattan}.  However, this is the extent of the similarity in approach between \cite{irving-rattan} and the current work: the details of our correspondence, and the trees themselves, differ from those in~\cite{irving-rattan}.

\begin{defn}
If a node in a ordered tree has any children, then it is a \emph{parent}.  
If a nonempty ordered tree contains at most one parent (the root), then it is a \emph{sapling}.
\end{defn}

\begin{ex}
Below are three examples of saplings.
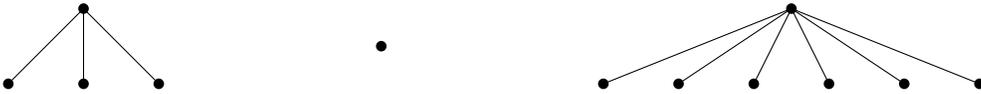
\begin{center}
\begin{tikzpicture}
\fill (0,0) circle (2pt);
\foreach \x in {-1,0,1} {\draw (0,0) -- (\x,-1); \fill (\x,-1) circle (2pt);}
\end{tikzpicture}
\hspace{1in}
\begin{tikzpicture}
\fill (0,0) circle (2pt);
\fill[white] (0,-.5) circle (2pt);
\end{tikzpicture}
\hspace{1in}
\begin{tikzpicture}
\fill (0,0) circle (2pt);
\foreach \x in {-2.5,-1.5,-.5,.5,1.5,2.5} {\draw (0,0) -- (\x,-1); \fill (\x,-1) circle (2pt);}
\end{tikzpicture}
\end{center}
\end{ex}

We will work with a set $\mc{T}_k(\pi)$ of trees, defined here and later, equivalently, in Definition~\ref{defn:tree hooks ornaments} (see Proposition~\ref{prop:tree is a bijection}).

\begin{defn}
Let $T_p$ be the sapling with $\ell_p$ nodes, where every leaf is labeled $p$.  For $i \in [m] \setminus \{p\}$, let $T_i$ be the sapling with $\ell_i+1$ nodes, where every leaf in the tree is labeled~$i$.
\end{defn}

We now describe the set $\mc{T}_k(\pi)$ of ordered trees specific to our work here.

\begin{defn}\label{defn:tree}
Consider the following iterative procedure.
\begin{itemize}
\item $T(0) = T_p$.
\item $T(j+1)$ is obtained from $T(j)$ by taking some $T_i$ that has not already been added, and inserting it into $T(j)$ by making the root of $T_i$ a new child of some parent node in $T(j)$, and giving this root the label~$i$.
\end{itemize}
Let $\mc{T}_{k}(\pi)$ consist of all possible $T(m-1)$ so obtained.
\end{defn}

\begin{ex}\label{ex:trees}
Consider $\pi = (142)(35)(6)$.  The following two trees are elements of $\mc{T}_{1}(\pi)$, where here $p=1$.
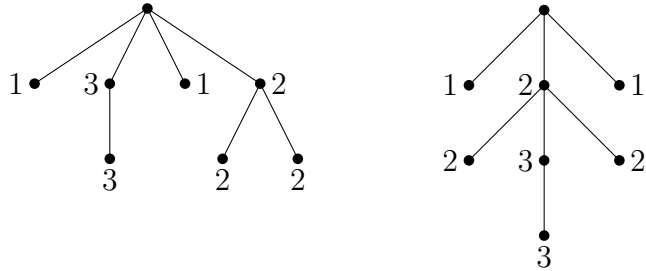
\begin{center}
\begin{tikzpicture}
\draw (0,0) -- (-1.5,-1) node[left] {$1$};
\draw (0,0) -- (-.5,-1) node[left] {$3$};
\draw (0,0) -- (.5,-1) node[right] {$1$};
\draw (0,0) -- (1.5,-1) node[right] {$2$};
\draw (-.5,-1) -- (-.5,-2) node[below] {$3$};
\draw (1.5,-1) -- (1,-2) node[below] {$2$};
\draw (1.5,-1) -- (2,-2) node[below] {$2$};
\fill (0,0) circle (2pt);
\foreach \x in {-1.5,-.5,.5,1.5} {\fill (\x,-1) circle (2pt);}
\foreach \x in {-.5,1,2} {\fill (\x,-2) circle (2pt);}
\draw (0,-3) node[below] {\phantom{!}};
\end{tikzpicture}
\hspace{.5in}
\begin{tikzpicture}
\draw (0,0) -- (-1,-1) node[left] {$1$};
\draw (0,0) -- (0,-3);
\draw (0,0) -- (1,-1) node[right] {$1$};
\draw (-1,-2) node[left] {$2$} -- (0,-1) -- (1,-2) node[right] {$2$};
\foreach \y in {0,-1,-2,-3} {\fill (0,\y) circle (2pt);}
\foreach \y in {-1,-2} {\foreach \x in {-1,1} {\fill (\x,\y) circle (2pt);};}
\draw (0,-1) node[left] {$2$};
\draw (0,-2) node[left] {$3$};
\draw (0,-3) node[below] {$3$};
\end{tikzpicture}
\end{center}
The following two trees are elements of $\mc{T}_{3}(\pi)$, where now $p=2$.
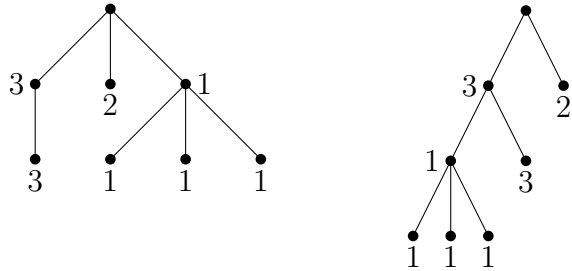
\begin{center}
\begin{tikzpicture}
\draw (0,0) -- (-1,-1) node[left] {$3$} -- (-1,-2) node[below] {$3$};
\draw (0,0) -- (0,-1) node[below] {$2$};
\draw (0,0) -- (1,-1) node[right] {$1$} -- (2,-2) node[below] {$1$};
\draw (0,-2) node[below] {$1$} -- (1,-1) -- (1,-2) node[below] {$1$};
\foreach \y in {0,-1,-2} {\fill (0,\y) circle (2pt);}
\foreach \y in {-1,-2} {\foreach \x in {-1,1} {\fill (\x,\y) circle (2pt);};}
\fill (2,-2) circle (2pt);
\draw (0,-3) node[below] {\phantom{!}};
\end{tikzpicture}
\hspace{.5in}
\begin{tikzpicture}
\draw (0,0) -- (-.5,-1) node[left] {$3$} -- (-1,-2) node[left] {$1$} -- (-1.5,-3) node[below] {$1$};
\draw (0,0) -- (.5,-1) node[below] {$2$};
\draw (-.5,-1) -- (0,-2) node[below] {$3$};
\draw (-1,-3) node[below] {$1$} -- (-1,-2) -- (-.5,-3) node[below] {$1$};
\fill (0,0) circle (2pt);
\foreach \x in {-.5,.5} {\fill (\x,-1) circle (2pt);}
\foreach \x in {-1,0} {\fill (\x,-2) circle (2pt);}
\foreach \x in {-1.5,-1,-.5} {\fill (\x,-3) circle (2pt);}
\end{tikzpicture}
\end{center}
\end{ex}

Each tree in the set $\mc{T}_{k}(\pi)$ corresponds to a word in $[m]^{n+m-2}$.

\begin{defn}\label{defn:tree to word}
Given $T \in \mc{T}_{k}(\pi)$, we obtain a word $\omega(T)$ by reading the labels of the non-root nodes in the order seen via a depth-first search.
\end{defn}

\begin{ex}
Continuing Example~\ref{ex:trees}, the first pair of trees, elements of $\mc{T}_{1}(\pi)$, map to the words $1331222$ and $1223321$ respectively, while the second pair of trees, elements of $\mc{T}_{3}(\pi)$, map to the words $3321111$ and $3111132$ respectively.
\end{ex}

The use of the letter $\omega$ to denote the maps in both Definitions~\ref{defn:decomp to word} and~\ref{defn:tree to word} is not coincidental.

\begin{prop}\label{prop:trees and words}
Given a word $\omega \in [m]^{n+m-2}$, we have $\omega = \omega(T)$ for some $T \in \mc{T}_{k}(\pi)$ if and only if $\omega = \omega({\delta})$ for some ${\delta} \in \star_{k}(\pi)$.
\end{prop}

\begin{proof}
This follows from Lemma~\ref{lem:allowable subwords} and Definition~\ref{defn:tree}.
\end{proof}

In other words, Proposition~\ref{prop:trees and words} describes a bijection between $\mc{T}_{k}(\pi)$ and $\mc{W}_{k}(\pi)$.  Although phrased in different terms, the set of valid words for $\star_{k}(\pi)$ is enumerated in~\cite{irving-rattan}:
$$|\mc{T}_{k}(\pi)| = |\mc{W}_{k}(\pi)| = (n+m-2)_{(m-2)}\ell_p.$$

\section{Bijective construction of minimal transitive star factorizations}\label{section:construction}

In this section, we use the characterization of $\star_k(\pi)$ of Proposition~\ref{prop:enough information} to give a map
$$\phi_{\pi,k} \ : \ \big\{[m-2] \hookrightarrow [n+m-2]\big\} \times [\ell_1] \times \cdots \times [\ell_m] \ \rightarrow \ \star_{k}(\pi),$$
where $\{X \hookrightarrow Y\}$ denotes the set of injections from $X$ into $Y$.  More precisely, we will obtain an element of $\mc{W}_{k}(\pi)$ and a set of cycle enclosures, which, by Proposition~\ref{prop:enough information}, define a unique element of $\star_{k}(\pi)$.  This element will be the image of~$\phi_{\pi,k}({x})$.

We will show that the map $\phi_{\pi,k}$ is a bijection for every $k$, thus obtaining a bijection from the set $\star_{k}(\pi)$ to the set $\star_{k'}(\pi)$: the composition of maps
$$\phi_{\pi,k'} \circ \left(\phi_{\pi,k}\right)^{-1}.$$

It is easiest to define the map $\phi_{\pi,k}$ via two preliminary operations.

\begin{defn}
Let $(\ell'_1, \ldots, \ell'_{m-1}) = (\ell_1,\ldots,\widehat{\ell_p},\ldots,\ell_m)$.  Fix
$$(c_1, \ldots, c_{m-1}) \in [\ell'_1] \times \cdots \times [\ell'_{m-1}],$$
and $i \in [m]\setminus \{p\}$.  If the $i$th cycle of $\pi$ in standard form is $(b_{i,1} \ b_{i,2} \ \cdots \ b_{i,\ell_i})$, then set
$$\ts{cycle}_{\pi,k}(c_1, \ldots, c_{m-1}) = \{b_{i,c_i} : i <p\} \cup \{b_{i,c_{i-1}} : i>p\}.$$
\end{defn}

We now compute $\ts{cycle}$ for our running example, as well as for a more complicated example that we will similarly examine throughout this section.

\begin{ex}\label{ex:cycle ex}
We have $\ts{cycle}_{(142)(35)(6),3}(2,1) = \{4,6\}$ because the second symbol in $(142)$ is $4$ and the first symbol in $(6)$ is~$6$.
\end{ex}

\begin{ex}\label{ex:cycle ex2}
We have $\ts{cycle}_{(18)(297)(3)(46)(5),9}(1,1,2,1) = \{1,3,5,6\}$ because the first symbol in $(18)$ is $1$, the first symbol in $(3)$ is $3$, the second symbol in $(46)$ is $6$, and the first symbol in $(5)$ is~$5$.
\end{ex}

The second operation, $\ts{tree}_{\pi,k}$, is rather more complex.  The idea is to take an injection $f:[m-2]\hookrightarrow[n+m-2]$ and a value $c \in [\ell_p]$, and to reinterpret them as a particular tree.  In this tree, all nodes except the root will have been labeled by values in~$[m]$.

\begin{defn}\label{defn:tree hooks ornaments}
Fix $(f,c) \in \{[m-2] \hookrightarrow [n+m-2]\} \times [\ell_p]$.  We now outline the procedure for producing the tree $\ts{tree}_{\pi,k}$.  
\begin{enumerate}
\item Label $[n+m-2] \setminus f([m-2])$, in increasing order, by ``$1$,'' $\ldots$, ``$2$,'' $\ldots$, ``$m,$'' $\ldots$, where each $i$ appears $\ell_i$ times.
\item Label the elements in the set $f([m-2])$ by the labels ``$f(i)$,'' as appropriate.
\item Change the $c$th occurrence of ``$p$'' to~``$f(0)$.''
\item Create $m$ \emph{factors} in $[n+m-2]$ by inserting $m-1$ bars: just after the rightmost ``$i$'' for each $i<p$ and just before the leftmost ``$i$'' for each~$i>p$.
\item Create a sapling $T_i$ with leaves labeled as in the $i$th factor.  Any leaves labeled by $f$ are \emph{hooks}, and the $T_i$ are \emph{ornaments}.
\item Let $T(0) = T_p$.
\item For $i \in [m-2]$, if the ornament $T_j$ containing ``$f(i)$'' has not yet been \emph{attached} in $T(i-1)$, then \emph{attach} it to $T(i-1)$ by identifying its root with the hook ``$f(i-1)$'' (called \emph{using} the hook), and label the identified node ``$j$,'' otherwise take no action; the resulting tree is~$T(i)$.  
\end{enumerate}
There remains at least 1 unattached ornament and at least 1 unused hook.  The attaching process and the fact that there were equally many ornaments and hooks at the start of the process means that there are the same number of unattached ornaments as unused hooks.  Let these be $T_{i_1}, \ldots, T_{i_r}$ and $f(h_1), \ldots, f(h_r)$, in increasing order of subscripts.  For each $j$, attach the root of $T_{i_j}$ to the hook ``$f(h_j)$'' as before.  The result is $\ts{tree}_{\pi,k}(f,c)$.
\end{defn}

We demonstrate Definition~\ref{defn:tree hooks ornaments} with two examples.

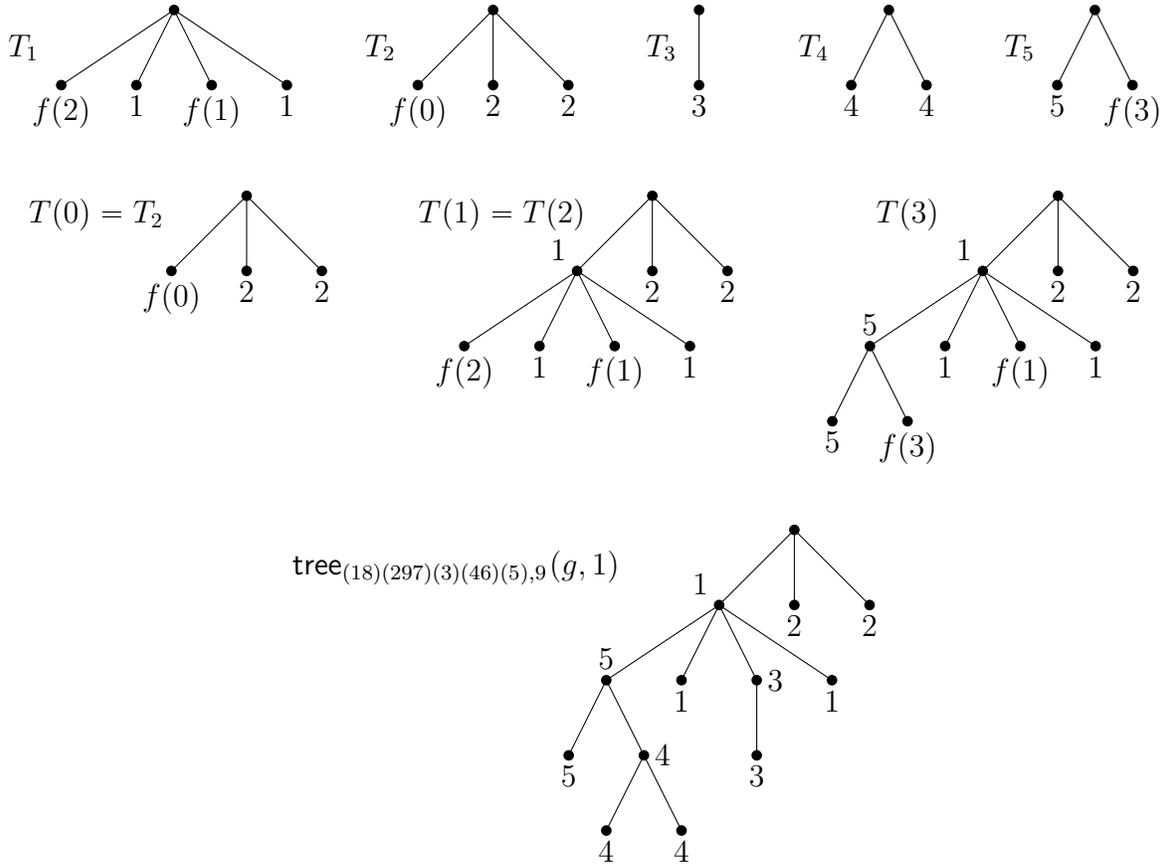
\begin{ex}\label{ex:tree ex}
Continuing Example~\ref{ex:cycle ex}, take $(142)(35)(6) \in \mf{S}_6$, where $n=6$ and $m=3$.  Let $k = 3$, so $p = 2$.  Let $f \in \{[1] \hookrightarrow [7]\}$ be defined by $f(1) = 3$, and let $c = 2$.   Definition~\ref{defn:tree hooks ornaments} produces the following work, where we represent the initial interval $[7]$ as a sequence of dots.
\begin{center}
\begin{tikzpicture}
\foreach \x in {1,2,3,4,5,6,7} {\fill (\x,0) circle (2pt);}
\foreach \x in {1,2,4} {\draw (\x,0) node[above] {$1$};}
\draw (3,0) node[above] {$f(1)$};
\draw (5,0) node[above] {$2$};
\draw (6,0) node[above] {$\xcancel{2}$};
\draw (6,.75) node {$f(0)$};
\draw (7,0) node[above] {$3$};
\draw[] (4.5,-.25) -- (4.5,.75);
\draw[] (6.5,-.25) -- (6.5,.75);
\end{tikzpicture}

\vspace{.25in}
\begin{tikzpicture}
\draw (.75,.25) node {$T_1$};
\draw (1,-.5) -- (2.5,.5) -- (4,-.5);
\draw (2,-.5) -- (2.5,.5) -- (3,-.5);
\foreach \x in {1,2,4} {\fill (\x,-.5) circle (2pt) node[below] {$1$};}
\fill (3,-.5) circle (2pt) node[below] {$f(1)$};
\fill (2.5,.5) circle (2pt);
\end{tikzpicture}
\hspace{.5in}
\begin{tikzpicture}
\draw (.75,.25) node {$T_2$};
\draw (1,-.5) -- (1.5,.5) -- (2,-.5);
\fill (1,-.5) circle (2pt) node[below] {$2$};
\fill (2,-.5) circle (2pt) node[below] {$f(0)$};
\fill (1.5,.5) circle (2pt);
\end{tikzpicture}
\hspace{.5in}
\begin{tikzpicture}
\draw (.5,.25) node {$T_3$};
\draw (1,-.5) -- (1,.5);
\fill (1,-.5) circle (2pt) node[below] {$3$};
\fill (1,.5) circle (2pt);
\draw (1,-.5) node[below] {\phantom{$f(1)$}};
\end{tikzpicture}

\vspace{.25in}

\begin{tikzpicture}
\draw (0,.25) node {$T(0) = T_2$};
\draw (1,-.5) -- (1.5,.5) -- (2,-.5);
\fill (1,-.5) circle (2pt) node[below] {$2$};
\fill (2,-.5) circle (2pt) node[below] {$f(0)$};
\fill (1.5,.5) circle (2pt);
\draw (1.5,-1.5) node[below] {\phantom{$f(1)$}};
\end{tikzpicture}
\hspace{.5in}
\begin{tikzpicture}
\draw (.5,.25) node {$T(1)$};
\draw (1,-.5) -- (1.5,.5) -- (2,-.5);
\fill (1,-.5) circle (2pt) node[below] {$2$};
\fill (2,-.5) circle (2pt) node[above right] {$1$};
\fill (1.5,.5) circle (2pt);
\draw (.5,-1.5) -- (2,-.5) -- (3.5,-1.5);
\draw (1.5,-1.5) -- (2,-.5) -- (2.5,-1.5);
\foreach \x in {.5,1.5,3.5} {\fill (\x,-1.5) circle (2pt) node[below] {$1$};}
\fill (2.5,-1.5) circle (2pt) node[below] {$f(1)$};
\end{tikzpicture}

\vspace{.25in}

\begin{tikzpicture}
\draw (-1.5,0) node {$\ts{tree}_{(142)(35)(6),3}(f, 2)$};
\draw (1,-.5) -- (1.5,.5) -- (2,-.5);
\fill (1,-.5) circle (2pt) node[below] {$2$};
\fill (2,-.5) circle (2pt) node[above right] {$1$};
\fill (1.5,.5) circle (2pt);
\draw (.5,-1.5) -- (2,-.5) -- (3.5,-1.5);
\draw (1.5,-1.5) -- (2,-.5) -- (2.5,-1.5);
\foreach \x in {.5,1.5,3.5} {\fill (\x,-1.5) circle (2pt) node[below] {$1$};}
\fill (2.5,-1.5) circle (2pt) node[right] {$3$};
\draw (2.5,-1.5) -- (2.5,-2.5);
\fill (2.5,-2.5) circle (2pt) node[below] {$3$};
\end{tikzpicture}

\end{center}
\end{ex}

\begin{ex}\label{ex:tree ex2}
Continuing Example~\ref{ex:cycle ex2}, take $(18)(297)(3)(46)(5) \in \mf{S}_9$, where $n = 9$ and $m=5$.  Let $k=9$, so $p=2$.  Let $g \in \{[3]\hookrightarrow[12]\}$ be defined by $g(1) = 3$, $g(2) = 1$, and $g(3) = 12$, and let $c = 1$.  Definition~\ref{defn:tree hooks ornaments} produces the following work, where we represent the initial interval $[12]$ as a sequence of dots.

\begin{center}
\begin{tikzpicture}
\foreach \x in {1,2,3,4,5,6,7,8,9,10,11,12} {\fill (\x,0) circle (2pt);}
\foreach \x in {2,4} {\draw (\x,0) node[above] {$1$};}
\foreach \x in {6,7} {\draw (\x,0) node[above] {$2$};}
\draw (8,0) node[above] {$3$};
\foreach \x in {9,10} {\draw (\x,0) node[above] {$4$};}
\draw (11,0) node[above] {$5$};
\draw (1,0) node[above] {$f(2)$};
\draw (3,0) node[above] {$f(1)$};
\draw (12,0) node[above] {$f(3)$};
\draw (5,0) node[above] {$\xcancel{2}$};
\draw (5,.75) node {$f(0)$};
\foreach \x in {4.5,7.5,8.5,10.5} {\draw[] (\x,-.25) -- (\x,.75);}
\end{tikzpicture}

\vspace{.25in}

\begin{tikzpicture}
\draw (.5,0) node {$T_1$};
\draw (1,-.5) -- (2.5,.5) -- (4,-.5);
\draw (2,-.5) -- (2.5,.5) -- (3,-.5);
\foreach \x in {2,4} {\fill (\x,-.5) circle (2pt) node[below] {$1$};}
\fill (3,-.5) circle (2pt) node[below] {$f(1)$};
\fill (1,-.5) circle (2pt) node[below] {$f(2)$};
\fill (2.5,.5) circle (2pt);
\end{tikzpicture}
\hspace{.15in}
\begin{tikzpicture}
\draw (.5,0) node {$T_2$};
\draw (1,-.5) -- (2,.5) -- (3,-.5); \draw (2,-.5) -- (2,.5);
\fill (1,-.5) circle (2pt) node[below] {$f(0)$};
\foreach \x in {2,3} {\fill (\x,-.5) circle (2pt) node[below] {$2$};}
\fill(2,.5) circle (2pt);
\end{tikzpicture}
\hspace{.15in}
\begin{tikzpicture}
\draw (.5,0) node {$T_3$};
\draw (1,-.5) -- (1,.5);
\fill (1,-.5) circle (2pt) node[below] {$3$};
\fill (1,.5) circle (2pt);
\draw (1,-.5) node[below] {\phantom{$f(3)$}};
\end{tikzpicture}
\hspace{.15in}
\begin{tikzpicture}
\draw (.5,0) node {$T_4$};
\draw (1,-.5) -- (1.5,.5) -- (2,-.5);
\foreach \x in {1,2} {\fill (\x,-.5) circle (2pt) node[below] {$4$};}
\fill (1.5,.5) circle (2pt);
\draw (1,-.5) node[below] {\phantom{$f(3)$}};
\end{tikzpicture}
\hspace{.15in}
\begin{tikzpicture}
\draw (.5,0) node {$T_5$};
\draw (1,-.5) -- (1.5,.5) -- (2,-.5);
\fill (1,-.5) circle (2pt) node[below] {$5$};
\fill (2,-.5) circle (2pt) node[below] {$f(3)$};
\fill (1.5,.5) circle (2pt);
\end{tikzpicture}

\vspace{.25in}

\begin{tikzpicture}
\draw (0,.25) node {$T(0) = T_2$};
\fill(2,.5) circle (2pt);
\draw (1,-.5) -- (2,.5) -- (3,-.5); \draw (2,-.5) -- (2,.5);
\fill (1,-.5) circle (2pt) node[below] {$f(0)$};
\foreach \x in {2,3} {\fill (\x,-.5) circle (2pt) node[below] {$2$};}
\draw (0,-2.5) node[below] {\phantom{$f(3)$}};
\end{tikzpicture}
\hspace{.25in}
\begin{tikzpicture}
\draw (0,.25) node {$T(1)=T(2)$};
\fill(2,.5) circle (2pt);
\draw (1,-.5) -- (2,.5) -- (3,-.5); \draw (2,-.5) -- (2,.5);
\fill (1,-.5) circle (2pt) node[above left] {$1$};
\foreach \x in {2,3} {\fill (\x,-.5) circle (2pt) node[below] {$2$};}
\draw (.5,-1.5) -- (1,-.5) -- (2.5,-1.5);
\draw (-.5,-1.5) -- (1,-.5) -- (1.5,-1.5);
\foreach \x in {.5,2.5} {\fill (\x,-1.5) circle (2pt) node[below] {$1$};}
\fill (1.5,-1.5) circle (2pt) node[below] {$f(1)$};
\fill (-.5,-1.5) circle (2pt) node[below] {$f(2)$};
\draw (0,-2.5) node[below] {\phantom{$f(3)$}};
\end{tikzpicture}
\hspace{.25in}
\begin{tikzpicture}
\draw (0,.25) node {$T(3)$};
\fill(2,.5) circle (2pt);
\draw (1,-.5) -- (2,.5) -- (3,-.5); \draw (2,-.5) -- (2,.5);
\fill (1,-.5) circle (2pt) node[above left] {$1$};
\foreach \x in {2,3} {\fill (\x,-.5) circle (2pt) node[below] {$2$};}
\draw (.5,-1.5) -- (1,-.5) -- (2.5,-1.5);
\draw (-.5,-1.5) -- (1,-.5) -- (1.5,-1.5);
\foreach \x in {.5,2.5} {\fill (\x,-1.5) circle (2pt) node[below] {$1$};}
\fill (1.5,-1.5) circle (2pt) node[below] {$f(1)$};
\fill (-.5,-1.5) circle (2pt) node[above] {$5$};
\draw (-1,-2.5) -- (-.5,-1.5) -- (0,-2.5);
\fill(-1,-2.5) circle (2pt) node[below] {$5$};
\fill(0,-2.5) circle (2pt) node[below] {$f(3)$};
\end{tikzpicture}

\vspace{.25in}

\begin{tikzpicture}
\draw (-2.5,0) node {$\ts{tree}_{(18)(297)(3)(46)(5),9}(g, 1)$};
\fill(2,.5) circle (2pt);
\draw (1,-.5) -- (2,.5) -- (3,-.5); \draw (2,-.5) -- (2,.5);
\fill (1,-.5) circle (2pt) node[above left] {$1$};
\foreach \x in {2,3} {\fill (\x,-.5) circle (2pt) node[below] {$2$};}
\draw (.5,-1.5) -- (1,-.5) -- (2.5,-1.5);
\draw (-.5,-1.5) -- (1,-.5) -- (1.5,-1.5);
\foreach \x in {.5,2.5} {\fill (\x,-1.5) circle (2pt) node[below] {$1$};}
\fill (1.5,-1.5) circle (2pt) node[right] {$3$};
\fill (-.5,-1.5) circle (2pt) node[above] {$5$};
\draw (-1,-2.5) -- (-.5,-1.5) -- (0,-2.5);
\fill(-1,-2.5) circle (2pt) node[below] {$5$};
\fill(0,-2.5) circle (2pt) node[right] {$4$};
\draw (1.5,-1.5) -- (1.5,-2.5);
\fill(1.5,-2.5) circle (2pt) node[below] {$3$};
\draw (-.5,-3.5) -- (0,-2.5) -- (.5,-3.5);
\foreach \x in {-.5,.5} {\fill(\x,-3.5) circle (2pt) node[below] {$4$};}
\end{tikzpicture}

\end{center}
\end{ex}

\begin{lem}
For $f \in \{[m-2] \hookrightarrow [n+m-2]\}$ and $c \in [\ell_p]$, we have $\ts{tree}_{\pi,k}(f,c) \in \mc{T}_{k}(\pi)$.
\end{lem}

\begin{proof}
We first must show that $\ts{tree}_{\pi,k}(f,c)$ is a tree.  What needs to be shown is that the last step, where unused hooks and unattached ornaments are identified, produces a tree.  At that stage, an ornament has not yet been attached to the tree if and only if none of its children are hooks.  Thus any unused hooks must appear in $T(m-2)$, so attaching hooks and ornaments as described does not create any cycles.

That this tree is an element of $\mc{T}_{k}(\pi)$ follows from Definition~\ref{defn:tree hooks ornaments}.
\end{proof}

\begin{prop}\label{prop:tree is a bijection}
The map $\ts{tree}_{\pi,k}$ is a bijection.
\end{prop}

\begin{proof}
We will show that the map is reversible.

Fix $T \in \mc{T}_{k}(\pi)$.  This tree can be decomposed into $m$ saplings based on the labels of the nodes.  Thus each sapling is associated with a value in $[m]$.  Call this the \emph{rank} of the sapling.  Let $S_i$ be the sapling of rank $i$.  (Note that the $m$ saplings can equivalently be identified by taking each parent together with all of its children that are leaves.)

If, in $T$, the root of the sapling $S$ is the parent of the root of the sapling $S'$, then $S$ \emph{shelters} $S'$ and write $S \succ S'$.  If a sapling shelters nothing, then it is \emph{free}.  Every sapling is sheltered by exactly one other sapling, except for $S_p$, which is sheltered by nothing.

Make maximal sequences $S_{j,1} \succ S_{j,2} \succ \ldots \succ S_{j,r_j}$; that is, $S_{j,1} = S_p$ for all $j$, and all $S_{j,r_j}$ are free.  Index the free saplings by $j$, so that their ranks are increasing with respect to $j$.  Let the roots of 
$$S_{1,2}, S_{1,3}, \ldots, S_{1,r_1}, S_{2,2}, S_{2,3}, \ldots, S_{2,r_2}, S_{3,2}, S_{3,3}, \ldots, S_{3,r_3}, \ldots$$
be relabeled ``$f(0)$,'' ``$f(1)$,'' $\ldots$, ``$f(m-2)$,'' respectively, with the provision that once a root has been labeled, the associated sapling is skipped when allocating the subsequent labels.  Let $T'$ be the resulting ordered tree after changing these $m-1$ labels in~$T$.

Now look in $T'$ at the $m$ saplings, and recall the ranks of each sapling as defined at the outset.  Working in increasing order of rank, write down the labels (as designated by $T'$) of all of the children of the root of each sapling from left to right (note that this includes non-leaf children).  Upon completion, we have written down $n-1+m-1 = n+m-2$ letters.  For $i\ge1$, if $f(i)$ appears in the $j$th position in this list, then set $f(i) = j$.  Set $c=f(0)$.  Then $T = \ts{tree}_{\pi,k}(f,c)$.
\end{proof}

We are now able to define the main bijection of this paper.

\begin{defn}\label{defn:phi}
Consider $(f, c_1, \ldots, c_m) \in \{[m-2] \hookrightarrow [n+m-2]\} \times [\ell_1] \times \cdots \times [\ell_m]$.  Let $\phi_{\pi,k}(f, c_1, \ldots, c_m)$ be the element of $\star_{k}(\pi)$ obtained, using Proposition~\ref{prop:enough information}, from the valid word $\omega(\ts{tree}_{\pi,k}(f,c_p))$ and the set of cycle enclosures $\ts{cycle}_{\pi,k}(c_1, \ldots, \widehat{c_p}, \ldots, c_m)$.
\end{defn}

To illustrate the map $\phi_{\pi,k}$, we continue the pair of examples studied throughout this section.

\begin{ex}\label{ex:phi k=3}
With the values determined in Examples~\ref{ex:cycle ex} and~\ref{ex:tree ex}, we have
$$\phi_{(142)(35)(6),3}(f,2,2,1) = (35)(34)(31)(32)(36)(36)(34) \in \star_{3}\big((142)(35)(6)\big).$$
\end{ex}

\begin{ex}\label{ex:phi 2}
With the values determined in Examples~\ref{ex:cycle ex2} and~\ref{ex:tree ex2}, we have
\begin{eqnarray*}
\phi_{(18)(297)(3)(46)(5),9}(g,1,1,1,2,1) &=& (91)(95)(95)(96)(94)(96)(98)(93)(93)(91)(92)(97)\\
&\in& \star_9 \big((18)(297)(3)(46)(5)\big).
\end{eqnarray*}
\end{ex}

\begin{thm}\label{thm:bijection}
For any $k \in [n]$, the map $\phi_{\pi,k}$ is a bijection.
\end{thm}

\begin{proof}
The sets $\{[m-2] \hookrightarrow [n+m-2]\} \times [\ell_1] \times \cdots \times [\ell_m]$ and $\star_{k}(\pi)$ have the same cardinality, by Theorem~\ref{thm:irving-rattan}.  The map $\ts{cycle}_{\pi,k}$ is certainly injective, and $\ts{tree}_{\pi,k}$ is a bijection, by Proposition~\ref{prop:tree is a bijection}.  The correspondence between pairs of valid words and cycle enclosures and elements of $\star_{k}(\pi)$ is a bijection, by Proposition~\ref{prop:enough information}.  

Therefore $\phi_{\pi,k}$ is a bijection.
\end{proof}

\section{Conclusion}\label{section:conclusion}

Theorem~\ref{thm:bijection} yields the desired combinatorial bijection between
$$\big\{[m-2] \hookrightarrow [n+m-2]\big\} \times [\ell_1] \times \cdots \times [\ell_m]$$
and $\star_k(\pi)$.  It also explains the symmetry of Theorem~\ref{thm:irving-rattan}.  That is, note that by Lemma~\ref{lem:allowable subwords}, no occurrence of $p$ in a valid word can sit between two occurrences of $j \neq p$.  Thus, if there are any other symbols appearing in a valid word (that is, if $m>1$), then the ``first'' of these must appear in one of the $\ell_p$ spots between or outside of the $\ell_p-1$ occurrences of $p$.  Choose $c_i \in [\ell_i]$ for each $i$.  Then $c_p$ determines where this ``first'' symbol appears relative to the $p$s in the valid word.  This information, together with the element of $\{[m-2] \hookrightarrow [n+m-2]\}$ yields the valid word via $\ts{tree}_{\pi,k}$.  For $i \neq p$, the $c_i$ determines the cycle enclosure of the $i$th cycle.  This explains that the role of the pivot $k$ affects only in what way the value $c_p \in [\ell_p]$ is interpreted by the bijection.

Thus, we have obtained a bijection
$$\phi_{\pi,k'} \circ \left(\phi_{\pi,k}\right)^{-1} \ : \ \star_k(\pi) \ \rightarrow \ \star_{k'}(\pi).$$

Finally, we demonstrate this bijection using the ongoing example of this article.

\begin{ex}
Let $k=3$, $\pi = (142)(35)(6)$, and $k'=1$.  Let us find the element of $\star_{k'}(\pi)$ corresponding to $(35)(34)(31)(32)(36)(36)(34) \in \star_k(\pi)$.  From Example~\ref{ex:phi k=3}, we know that
$$\left(\phi_{\pi,k}\right)^{-1} \big((35)(34)(31)(32)(36)(36)(34)\big) = (f:1\mapsto3,2,2,1).$$
Now we apply $\phi_{\pi,k'}$ to $(f:1\mapsto3,2,2,1)$, obtaining
$$\phi_{\pi,k'} \circ \left(\phi_{\pi,k}\right)^{-1} \big((35)(34)(31)(32)(36)(36)(34)\big) = (12)(15)(13)(15)(16)(16)(14) \in \star_{k'}(\pi).$$
\end{ex}

\end{document}